\documentclass{amsart} 
\title{ A criterion to specify the absence of Baire property }
\theoremstyle{plain}
\newtheorem{theorem}{Theorem}[section]

\newtheorem*{corollary*}{Corollary}

\theoremstyle{definition}

\newtheorem{example}[theorem]{Example}

\usepackage{tikz}

\author{Mehdi Pourbarat}

\address{Shahid Beheshti University, Department of Mathematics,Tehran, Iran}
\email{m-pourbarat@sbu.ac.ir}
\author{Neda Abbasi}

\address{Shahid Beheshti University, Department of Mathematics,Tehran, Iran}
\email{n\_abbasi@sbu.ac.ir}

\thanks{}
\date{\today}
\usepackage[pagebackref=true,colorlinks,linkcolor=blue,citecolor=magenta]{hyperref}
\usepackage{yhmath,graphicx,amsmath,amsfonts,amssymb,amsthm,color,amscd,mathtools} 
\usepackage[normalem]{ulem}

\begin{document}

\begin{abstract} Let $X$ be a topological space. Let $X_0 \subseteq X$ be a second countable subspace. Also, assume that $X$ is first countable at any point of $X_0$. Then we provide
some conditions under which we ensure that $X_0$ is not Baire. 

\end{abstract}
\maketitle Subject Classification: 37B55, 54E52. \\
\maketitle Keywords: Nonautonomous systems, topological transitivity, Baire space, Birkhoff theorem.
\maketitle

\maketitle
\section{Introduction }
A space $X$ is called Baire if the intersection of any sequence of dense open subsets of $X$ is dense in $X$.
Alternatively, this notation can be formulated in terms of second category sets.
The Baire category theory has numerous applications in Analysis and Topology.
Among these applications are, for instance, the open mapping, closed graph theorem and
the Banach-Steinhaus theorem in Functional Analysis \cite {{Aarts}, {Haworth}}.

 The aim of this paper is to introduce a trick that concludes the absence of Baire property for some  topological spaces using dynamical techniques and tools. Before
stating the main result, we establish some notations.

Let $(X,~\tau)$ be a topological space, $X_n$'s its subspaces and 
$$ x_{n+1}=f_n(x_n), ~ n \in \mathbb{N} \cup \{0\}, $$ where
$ f_n: {X_n} \to { X_{n+1}}$ are  continuous maps.
The family $\{f_n\}_{n=0}^{\infty}$ is called a nonautonomous discrete system \cite{{Shi}, {Shi-Chen}}.
For given $x_{0}\in X_{0}$, the orbit of $x_0$ is  defined as
$$ orb(x_{0}):=\big\{ x_{0}, ~ f_{0}(x_{0}), ~ f_{1}\circ f_{0}(x_{0}), ~\cdots , ~f_{n} \circ f_{n-1} \circ \cdots \circ f_{0}(x_{0}),~\cdots \big\},$$
and we  say that this orbit starts from the point $x_0$. The topological structure of  the orbit that starts from the point $x_0$ may be  complex.
Here, we study the points of $ X_{0}$ whose orbits always intersect around $ X_{0}$. They are formulated as follows:
$$O:= \big\{x \in X_0 | ~~ {\overline{orb(x)}}^{{X}}\cap X_0=X_0 \big\}.$$
The system $\{f_n\}_{n=0}^{\infty}$ is called topologically transitive on $ X_{0}$  if for any two non-empty open sets $U_{0}$ and $ V_{0}$ in $ X_{0}, $ there
exists $ n \in \mathbb{N}$ such that $U_{n}\cap V_{0}\neq \phi $, where $ U_{i+1}=f_{i}(U_i)$ for $ 0 \leq i \leq n-1$, in other word
$(f_{n-1} \circ f_{n-2} \circ \cdots \circ f_{1} \circ f_{0})(U_{0}) \cap V_{0}\neq \phi$ \cite {Shi-Chen}.

Our main theorem is as follows:
\begin{theorem} \label{Theorem }
 Let $X$ be a topological space. Let $X_0$ be a second countable subspace of $ X$ and let $X$ be first countable at any point of $X_0$. Also, suppose that
 the system $\{f_n\}_{n=0}^{\infty}$ is topologically transitive on $X_{0}$ and  $ \overline{O}\neq X_0$. Then $X_0$ can  not  be a Baire subspace.
\end{theorem}
Note that, if $X$ is a metric space, $X_n=X$, and $f_n=f $ for each $n$, then Theorem \ref{Theorem }  will be obtained as  a direct result of Birkhoff transitivity theorem. This fact was   our   motivation in writing the paper.
\section{Proof }

So as $X_{0}$ is  a second countable subspace and $X$ first countable at any point of $X_0$, it is easy to show that
there exists a collection $\{U_m\}_{m \in N}$ of open sets in $X$ such that
\begin{itemize}

\item [$i$)] $ U_m \cap D_0 \neq\phi$,

\item [$ii$)]  the family $\{U_m \cap D_0\}_{m \in \mathbb{N}}$ is a basis for $D_0$,

\item [$iii$)]for each $x_0\in D_0$, the family  $\{U_m\}_{m \in {\mathbb{N}}}$ is a local basis for  $x_0$ in $X$.
\end{itemize}
We  claim that 
$$O=\bigcap_{m=1}^{\infty}\bigcup_{n=1}^{\infty}{f_{n-1}}
 \circ {f_{n-2}} \circ \cdots \circ {f_{1}}  \circ {f_{0}}^{-1}(U_{m}). \eqno{(2.1)}$$
To prove the claim, put $O^*:=\bigcap_{m=1}^{\infty}\bigcup_{n=1}^{\infty}{f_{n-1}\circ f_{n-2} \circ \cdots \circ f_{1}  \circ f_{0}}^{-1}(U_{m}).$
Firstly, we show that $O \subseteq O^*$. Suppose otherwise, there is $ x \in O $ such that $x \notin O^*$. So as $x \notin O^*$, 
there exists $m \in \mathbb{N}$ such that for each $n \in \mathbb{N}$ we have 
$$ {(f_{n-1}\circ f_{n-2}\circ \cdots \circ f_{1} \circ f_{0})}(x) \notin {U_{m}}.$$
Hence, $orb(x) \cap U_m=\phi$. Since $U_m\cap D_0 \neq \phi$, there exists an element $z\in U_m\cap X_0$, such that $z\notin {\overline{orb(x)}}^X.$ 
But $ z \in X_0$ and so $ {\overline{orb(x)}}^X\cap X_0 \neq D_0$. It is concluded that $x \notin O$ which contradicts the choice of $x$.
Now, it is shown that $O^* \subseteq O$. Let $x\in O^*$ but $x \notin O$. So as $x\in O^*$, concluded
for each $m \in \mathbb{N}$, there exists $n \in \mathbb{N}$ such that 
$ {(f_{n-1}\circ f_{n-2} \circ \cdots \circ f_{1} \circ f_{0})}(x) \in U_m.$
Thus, $orb(x) \cap U_m \neq \phi$. Moreover, the relation $x \notin O$ indicates that there exists $ z\in X_0$ such that $z \notin \overline{orb(x)}^X.$ Consequently,
there exists $U_k$ containing $z$, such that $U_k \cap orb(x)= \phi$ that this contradicts with $ orb(x) \cap U_m \neq \phi$, for each $m\in\mathbb{N}$. 
\\ By continuity of $ f_n: {X_n} \rightarrow { X_{n+1}},$ 
each set $\bigcup_{n=1}^{\infty}{(f_{n-1} \circ f_{n-2} \circ \cdots  \circ f_{1} \circ f_{0})}^{-1}(U_{m})$ is open and because of transitivity, these open sets are dense in $X_{0}$. If $X_{0}$ be  a Baire space,  then (2.1) implies that $O$ is a dense $ G_{\delta}$-set.  This  is a contradict with $ \overline{O}\neq X_0$. Thus $X_0$ is not  a Baire subspace, and the proof of the Theorem \ref{Theorem } is complete. 

\section{Example }
\begin{example}
 Consider  $X=\mathbb H(\mathbb C)
=\big\{f:\mathbb C \rightarrow \mathbb C|~ f ~ is ~holomorphic \big\}$ endowed with the metric $d(f,g)=\displaystyle\sum_ {n=1}^{\infty}\frac {1}{2^n} min \big(~1,~p_n(f-g)\big),$
with $p_n(h) ={sup}_{|z| \leq  n } |h(z)|$. 
Then $X$ is a separable Banach space and besides that the differentiation operator $D:\mathbb H(\mathbb C) \rightarrow \mathbb H(\mathbb C)$ with $ D(f)=f^\prime$ is continuous \cite {Grosse-Erdmann}. Moreover, the space $\mathbb H(\mathbb C)$ is Baire and if we  consider the  dynamical system $ D:\mathbb H(\mathbb C ) \rightarrow \mathbb H(\mathbb C),$ then Birkhoff theorem guarantees the existence of  functions  that their orbit is dense in $\mathbb{ H} (\mathbb{ C})$.
  Now,  assume that 
$$ X_0=\big\{\sum_{i=0}^{N} a_iz^i+\alpha g(z)\big |~ a_i , \alpha \in \mathbb{C} \big\}.$$
Then the subspace $ X_0$ is not  Baire.  To see this,  take $\{\alpha_n\}_{n=0}^{\infty}$ be a subsequence
 with $\alpha_0=0$ in this way that $D^{\alpha_n}(g)$ is convergent. We consider nonautonomous discrete system $\{f_{n}\}_{n=0}^{\infty}$ with $f_n=D^{\alpha_{n+1}-\alpha_n}$  where 
$ X_{n}=\big\{\sum_{i=0}^{N} a_iz^i+\alpha  g^{(\alpha_n)}(z)\big | a_i , \alpha \in \mathbb{C} \big\}$. By planning the arguments  similar to what employed in  the proof of Example 2.21 in \cite{Grosse-Erdmann}, we observe that the 
system $\{f_n\}_{n=0}^{\infty}$ is topologically transitive. Now the assertion obtains  by using Theorem \ref{Theorem } since the set $ O$ is empty.

\end{example}

\end{document}